\documentclass[12pt]{amsart}

\usepackage{amssymb}
\usepackage{amstext}
\usepackage{amsmath}
\usepackage{amscd}
\usepackage{latexsym}
\usepackage{amsfonts}

\newtheorem{thm}{Theorem}[section]

\newtheorem{cor}[thm]{Corollary}

\newtheorem{lem}[thm]{Lemma}

\newtheorem{ex}[thm]{Example}

\newcommand{\gp}{\mathfrak{p}}

\newcommand{\gm}{\mathfrak{m}}

\newcommand{\ga}{\mathfrak{a}}
\newcommand{\gb}{\mathfrak{b}}

\newcommand{\mult}[2]{{\mathrm e}_{#1}(#2)} 
\newcommand{\nmult}[2]{{\overline{\mathrm e}}_{#1}(#2)} 

\newcommand{\length}[2]{\ell_{#1}(\,#2\,)}

\newcommand{\height}[2]{{\rm ht}_{#1}\,#2}
\newcommand{\dep}[2]{{\rm depth}_{#1}\,#2}

\newcommand{\rn}[2]{{\rm r}_{#1}(#2)} 
\newcommand{\nog}[2]{\mu_{#1}(#2)} 
\newcommand{\rrc}[1]{\widetilde{#1}} 
\newcommand{\ic}[1]{\overline{#1}} 

\begin{document}

\setlength{\baselineskip}{16pt}

\title{An upper bound on the reduction number of an ideal}

\author{Yayoi Kinoshita}
\address{Yayoi Kinoshita,
Department of Mathematics and Informatics,
Graduate School of Science, Chiba University,
1-33 Yayoi-cho, Inage-ku, Chiba-shi, 263-8522 Japan}
\email{}

\author{Koji Nishida}
\address{Koji Nishida (corresponding author),
Department of Mathematics and Informatics,
Graduate School of Science, Chiba University,
1-33 Yayoi-cho, Inage-ku, Chiba-shi, 263-8522 Japan}
\email{nishida@math.s.chiba-u.ac.jp}

\author{Kensuke Sakata}
\address{Kensuke Sakata,
Department of Mathematics and Informatics,
Graduate School of Science, Chiba University,
1-33 Yayoi-cho, Inage-ku, Chiba-shi, 263-8522 Japan}
\email{}

\author{Ryuta Shinya}
\address{Ryuta Shinya,
Department of Mathematics and Informatics,
Graduate School of Science, Chiba University,
1-33 Yayoi-cho, Inage-ku, Chiba-shi, 263-8522 Japan}
\email{}

\thanks{
\endgraf
}

\begin{abstract}
Let $A$ be a commutative ring and $I$ an ideal of $A$
with a reduction $Q$.
In this paper we give an upper bound on the reduction number
of $I$ with respect to $Q$,
when a suitable family of ideals in $A$ is given.
As a corollary it follows that if some ideal $J$
containing $I$ satisfies $J^2 = QJ$,
then $I^{v + 2} = QI^{v + 1}$, where $v$ denotes
the number of generators of $J / I$
as an $A$-module.
\end{abstract}

\maketitle

\section{Introduction}
Let $Q$, $I$ and $J$ be ideals of a commutative ring $A$
such that $Q \subseteq I \subseteq J$.
As is noted in \cite[2.6]{GNO},
if $J / I$ is cyclic as an $A$-module and $J^2 = QJ$,
then we have $I^3 = QI^2$.
The purpose of this paper is to generalize this fact.
We will show that if $J / I$ is generated by $v$ elements
as an $A$-module and $J^2 = QJ$,
then $I^{v + 2} = QI^{v + 1}$.
We get this result as a corollary of the following theorem,
which generalizes Rossi's assertion stated in the proof
of \cite[1.3]{R}.

\begin{thm}\label{1a}
Let $A$ be a commutative ring and $\{ F_n \}_{n \geq 0}$
a family of ideals in $A$ such that
$F_0 = A$, $IF_n \subseteq F_{n + 1}$ for any $n \geq 0$,
and $I^{k + 1} \subseteq QF_k + \ga F_{k + 1}$
for some $k \geq 0$ and an ideal $\ga$ in $A$.
Suppose that $F_n / (QF_{n - 1} + I^n)$ is
generated by $v_n$ elements for any $n \geq 0$
and $v_n = 0$ for $n \gg 0$.
We put $v = \sum_{n \geq 0} v_n$.
Then we have
\[
I^{v + k + 1} = QI^{v + k} + \ga I^{v + k + 1}\,.
\]
\end{thm}

If a family $\{ F_n \}_{n \geq 0}$ of ideals
in $A$ satisfies all of the conditions required in \ref{1a}
in the case where $\ga = (0)$,
we have $F_n = QF_{n - 1}$ for $n \gg 0$.
As a typical example of such $\{ F_n \}_{n \geq 0}$,
we find $\{ \rrc{I^n} \}_{n \geq 0}$
when $I$ contains a non-zerodivisor,
where $\widetilde{I^n}$ denotes the Ratliff-Rush closure of $I^n$
(cf. \cite{RR}).
If $A$ is an analytically unramified local ring,
then $\{ \overline{I^n} \}_{n \geq 0}$ is also an important example,
where $\ic{I^n}$ denotes the integral closure of $I^n$.
It is obvious that $\{ J^n \}_{n \geq 0}$ always satisfies
the required condition on $\{ F_n \}_{n \geq 0}$ for any
ideal $J$ with $I \subseteq J \subseteq \ic{I}$.

We prove \ref{1a} following Rossi's argument
in the proof of \cite[1.3]{R}.
However we do not assume that $A / I$ has finite length.
And furthermore we can deduce the following corollary which
gives an upper bound on the reduction number
$\rn{Q}{I}$ of $I$ with respect to $Q$
using numbers of gerators of certain $A$-modules.

\begin{cor}\label{1b}
Let $(A, \gm)$ be a Noetherian local ring and
$\{ F_n \}_{n \geq 0}$ a family of ideals
in $A$ such that $F_0 = A$, $IF_n \subseteq F_{n + 1}$
for any $n \geq 0$, and
$I^{k + 1} \subseteq QF_k + \gm F_{k + 1}$ for some $k \geq 0$.
Then we have
\begin{eqnarray*}
\rn{Q}{I} & \leq & k + \sum_{n \geq 1}\, \nog{A}{F_n / (QF_{n - 1} + I^n)} \\
 & \leq & 1 + \nog{A}{F_1 / I} + \sum_{n \geq 2}\, \nog{A}{F_n / QF_{n - 1}}\,.
\end{eqnarray*}
\end{cor}

Throughout this paper $A$ denotes a commutative ring.
We do not assume that $A$ is Noetherian unless otherwise specified.
Furthermore $I$ and $Q$ denote ideals of $A$ such that $Q \subseteq I$.
We set $\rn{Q}{I} = \inf\{ n \geq 0 \mid I^{n + 1} = QI^n \}$.
Of course, $\rn{Q}{I} = \infty$ if
$Q$ is not a reduction of $I$.
For a finitely generated $A$-module $M$, we denote by $\nog{A}{M}$
the minimal number of generators of $M$.
If $(A, \gm )$ is a Noetherian local ring
and $M$ is annihilated by some power of $\gm$,
the length of $M$ is denoted by $\length{A}{M}$.

\vspace{1em}

\section{Proof of Theorem \ref{1a}}
In order to prove \ref{1a} we prepare
the following lemma, which generalizes \cite[2.3]{Hc}.

\begin{lem}\label{2a}
Let $I_1, I_2, \dots\,, I_N$ be finite number of ideals of $A$.
For any $1 \leq n \leq N$,
we assume that $I_n$ is generated by $v_n$ elements and
\[
I \cdot I_n \subseteq I^{n+1} + \sum_{\ell = 1}^N Q^{n + 1 - \ell}I_\ell\,.
\]
Let $v := v_1 + v_2 + \cdots + v_N > 0$.
Then, for any $v$ elements $a_1, a_2, \dots\,, a_v$ in $I$,
there exists $\sigma \in QI^{v -1}$ such that
\[
a_1a_2 \cdots a_v - \sigma \in \bigcap_{n = 1}^N\,[I^{n + v} : I_n]\,.
\]
\end{lem}

\noindent
{\it Proof.}\hspace{0.5ex}
We put $w_0 = 0$ and $w_n = v_1 + \cdots + v_n$
for $1 \leq n \leq N$.
Then $0 = w_0 \leq w_1 \leq w_2 \leq \cdots \leq w_N = v$.
Hence, if $1 \leq i \leq v$,
we have $w_{n - 1} < i \leq w_n$ for some $1 \leq n \leq N$,
and we denote this number $n$ by $n_i$.
Now we choose elements $x_1, x_2, \dots\,, x_v$ of $A$ so that
$I_n$ is generated by $\{ x_i \mid w_{n-1} < i \leq w_n \}$
for any $1 \leq n \leq N$ with $v_n \neq 0$.
Then $x_i \in I_{n_i}$ and
\[
a_ix_i \in I \cdot I_{n_i} \subseteq
I^{n_i + 1} + \sum_{\ell = 1}^N\,Q^{n_i + 1 - \ell}I_\ell
\]
for any $1 \leq i \leq v$.
Hence there exists a family $\{ c_{ij} \}_{1 \leq i, j \leq v}$
of elements in $A$ such that
\[
a_ix_i \equiv \sum_{j = 1}^v\,c_{ij}x_j \hspace{1ex} \mbox{mod} \hspace{1ex} I^{{n_i} + 1}
\hspace{2ex} \mbox{and} \hspace{2ex}
c_{ij} \in Q^{n_i + 1 - n_j}
\]
for any $1 \leq i, j \leq v$.
Let $R = A[It, t^{-1}]$ and $T = A[t, t^{-1}]$,
where $t$ is an indeterminate.
We regard $T / R$ as a graded $R$-module,
and for any $f \in T$ we denote by $\ic{f}$ the class of $f$ in $T / R$.
Then we have
\[
a_it \cdot \ic{x_it^{n_i}} = \sum_{j = 1}^v\,
c_{ij}t^{n_i - n_j + 1} \cdot \ic{x_jt^{n_j}}
\]
for any $1 \leq i, j \leq v$.
Here we put
\begin{eqnarray*}
b_{ij} & = & \left\{\begin{array}{ll}
a_i - c_{ii} & \mbox{if $i = j$} \\
-c_{ij}      & \mbox{if $i \neq j$}
\end{array}\right.
\,, \\
m_{ij} & =  & b_{ij}t^{n_i - n_j + 1} \in R\,, \hspace{1ex} \mbox{and} \\
e_i & = & \ic{x_it^{n_i}} \in T / R
\end{eqnarray*}
for any $1 \leq i, j \leq v$.
Let us consider the $v \times v$ matrix
$M = (\,m_{ij}\,)$ with entries in $R$.
Because we have
\[
M
\left(\begin{array}{c}
e_1 \\
e_2 \\
\vdots \\
e_v
\end{array}\right)
=
\left(\begin{array}{c}
0 \\
0 \\
\vdots \\
0
\end{array}\right) \,,
\]
it follows that $\Delta e_i = 0$ for any $1 \leq i \leq v$,
where $\Delta = \det M$.
Then we get
\begin{equation}
\Delta \cdot x_it^{n_i} \in R \tag{$\ast$}
\end{equation}
for any $1 \leq i \leq v$.
On the other hand, by the definition of determinant, we have
\[
\Delta = \sum_{(p_1, p_2, \dots, p_v) \in S_v}\,
{\rm sgn}(p_1, p_2, \dots, p_v) m_{1p_1}m_{2p_2} \cdots m_{vp_v}\,,
\]
where $S_v$ denotes the set of permutations of $1, 2, \dots, v$
and ${\rm sgn}(p_1, p_2, \dots, p_v)$ denotes the signature of
$(p_1, p_2, \dots, p_v) \in S_v$.
Because
\[
\deg (\prod_{i = 1}^v m_{ip_i})
= \sum_{i = 1}^v (n_i - n_{p_i} + 1)
= \sum_{i = 1}^v n_i - \sum_{i = 1}^v n_{p_i} + v
= v\,,
\]
we have $\prod_{i = 1}^v m_{ip_i} = (\prod_{i = 1}^v b_{ip_i})t^v$.
Therefore $\Delta = \delta t^v$, where $\delta$ denotes the determinant of
the $v \times v$ matrix $(\,b_{ij}\,)$ with entries in $A$.
Hence, by $(\ast)$ we have $\delta x_i \in I^{v + n_i}$
for any $1 \leq i \leq v$.
This means $\delta I_n \subseteq I^{v + n}$
for any $1 \leq n \leq N$,
and so $\delta \in \bigcap_{n = 1}^N [I^{v + n} : I_n]$.
If $(p_1, p_2, \dots\,, p_v) \neq (1, 2, \dots\,, v)$,
then $j > p_j$ for some $1 \leq j \leq v$, which means
$n_j \geq n_{p_j}$, and so
$b_{jp_j} = - c_{jp_j} \in Q^{n_j - n_{p_j} + 1} \subseteq QI^{n_j - n_{p_j}}$.
As a consequence, if $(1, 2, \dots, v) \neq (p_1, p_2, \dots, p_v) \in S_v$,
we get 
\[
\prod_{i = 1}^v b_{ip_i}
= b_{jp_j} \cdot \prod_{i \neq j} b_{ip_i}
= QI^{n_j - n_{p_j}} \cdot \prod_{i \neq j} I^{n_i - n_{p_i} + 1}
\subseteq Q \cdot I^{n_j - n_{p_j} + \sum_{i \neq j}(n_i - n_{p_i} + 1)}
= QI^{v - 1} \,.
\]
Furthermore, as $a_i \in I$ and $c_{ii} \in Q$
for any $1 \leq i \leq v$, we have
\[
\prod_{i = 1}^v b_{ii} = \prod_{i = 1}^v (a_i - c_{ii}) =
a_1a_2 \cdots a_v - d
\]
for some $d \in QI^{v - 1}$.
Therefore, there exists $\sigma \in QI^{v - 1}$ such that
$\delta = a_1a_2 \cdots a_v - \sigma$,
and the proof is complete.
 
\vspace{1.5em}
\noindent
{\it Proof of Theorem 1.1.} \hspace{0.5ex}
If $v = 0$,
then we have $F_n = I^n$ for any $n \geq 0$,
and so $I^{k + 1} \subseteq QF_k + \ga F_{k + 1} =
QI^k + \ga I^{k + 1} \subseteq I^{k + 1}$,
which means $I^{k + 1} = QI^k + \ga I^{k + 1}$.
Hence we may assume $v > 0$.
For any $n \geq 0$, let us take an ideal $I_n$
generated by $v_n$ elements so that
$F_n = QF_{n - 1} + I^n + I_n$.
We can easily show that
\begin{equation}
F_n = I^n + \sum_{\ell = 0}^n Q^{n - \ell}I_\ell \tag{\#}
\end{equation}
for any $n \geq 0$ by induction on $n$.
Now we choose an integer $N$ so that $N > k$ and
$I_n = 0$ for any $n > N$.
Then by $(\#)$ it follows that
\[
I \cdot I_n \subseteq F_{n + 1} =
I^{n + 1} + \sum_{\ell = 0}^N Q^{n + 1 - \ell}I_\ell
\]
for any $0 \leq n \leq N$.
Let $a_1, a_2, \dots, a_v$ be any elements of $I$.
Then, by \ref{2a} there exists $\sigma \in QI^{v - 1}$ such that
\[
a_1a_2 \cdots a_v - \sigma \in \bigcap_{n = 0}^N\, [I^{n + v} : I_n]\,.
\]
We put $\xi = a_1a_2 \cdots a_v - \sigma$.
Then by $(\#)$ we get
\[
\xi F_n = \xi I^n + \sum_{\ell = 0}^n Q^{n - \ell} \cdot \xi I_\ell
\subseteq I^v \cdot I^n + \sum_{\ell = 0}^n Q^{n - \ell} \cdot I^{\ell + v}
\subseteq I^{v + n}
\]
for any $0 \leq n \leq N$.
Now the assumption that
$I^{k + 1} \subseteq QF_k + \ga F_{k + 1}$ implies
\[
\xi I^{k + 1} \subseteq Q \cdot \xi F_k + \ga \cdot \xi F_{k + 1}
\subseteq Q \cdot I^{v + k} + \ga \cdot I^{v + k + 1}\,.
\]
Therefore we get
\[
a_1a_2 \cdots a_v \cdot I^{k + 1} = (\xi + \sigma)I^{k + 1} \subseteq
QI^{v + k} + \ga I^{v + k + 1}\,.
\]
Then, as the elements $a_1, a_2, \dots, a_v$ are chosen arbitrarily from $I$,
it follows that $I^v \cdot I^{k + 1} \subseteq
QI^{v + k} + \ga I^{v + k + 1}
\subseteq I^{v + k + 1}$.
Thus we get $I^{v + k + 1} = QI^{v + k} + \ga I^{v + k + 1}$.

\vspace{1.2em}
\noindent
{\it Proof of Corollary 1.2.} \hspace{0.5ex}
We put $v = \sum_{n \geq 1} \nog{A}{F_n / (QF_{n - 1} + I^n)}$.
We may assume $v < \infty$.
Then, setting $\ga = \gm$ in \ref{1a},
it follows that
$I^{v + k + 1} = QI^{v + k} + \gm I^{v + k + 1}$.
Hence we get $I^{v + k + 1} = QI^{v + k}$
by Nakayama's lemma,
and so $\rn{Q}{I} \leq v + k$.
In order to prove the second inequality,
we choose $k$ as small as possible.
If $k \leq 1$, we have
\[
\rn{Q}{I} \leq k + v
\leq 1 + \nog{A}{F_1 / I} + \sum_{n \geq 2} \nog{A}{F_n / QF_{n - 1}}\,.
\]
So, we assume $k \geq 2$ in the rest of this proof.
In this case we have
\begin{equation}
\rn{Q}{I} \leq k + \nog{A}{F_1 / I} +
\sum_{n = 2}^k \nog{A}{F_n / (QF_{n - 1} + I^n)} +
\sum_{n \geq k + 1} \nog{A}{F_n / QF_{n - 1}}\,. \tag{$\natural$}
\end{equation}
If $2 \leq n \leq k$,
then $I^n \not\subseteq QF_{n - 1} + \gm F_n$,
and so the canonical surjection
\[
F_n / (QF_{n - 1} + \gm F_n) \longrightarrow
F_n / (QF_{n - 1} + I^n + \gm F_n)
\]
is not injective, which means
\[
\nog{A}{F_n / QF_{n - 1} + I^n} \leq \nog{A}{F_n / QF_{n - 1}} - 1\,.
\]
Thus we get
\[
\sum_{n = 2}^k \nog{A}{F_n / QF_{n - 1} + I^n} \leq
\{\sum_{n = 2}^k \nog{A}{F_n / QF_{n - 1}}\} - (k - 1)\,.
\]
Therefore the required inequality follows from $(\natural)$.

\vspace{1em}

\section{Corollaries}
In this section we collect some results deduced
from \ref{1a} and \ref{1b}.

\begin{cor}\label{3a}
Let $J$ be an ideal of $A$ such that
$J \supseteq I$ and $J^2 = QJ$.
If $J / I$ is finitely generated as an
$A$-module,
then $\rn{Q}{I} \leq \nog{A}{J / I} + 1$.
\end{cor}

\noindent
{\it Proof.}\hspace{0.5ex}
We apply \ref{1a} setting $F_n = J^n$ for any $n \geq 0$
and $\ga = (0)$.
Because $I^2 \subseteq J^2 = QJ$,
we may put $k = 1$, and hence we get
$I^{v + 2} = QI^{v + 1}$,
where $v = \nog{A}{J / I}$.
Then $\rn{Q}{I} \leq v + 1$.

\begin{cor}\label{3b}
Let $(A, \gm)$ be a two-dimensional regular local ring
{\rm (}or, more generally, a two-dimensional pseudo-rational
local ring{\rm )} such that $A / \gm$ is infinite.
If $I$ is an $\gm$-primary ideal with a minimal reduction $Q$,
then $\rn{Q}{I} \leq \nog{A}{\ic{I} / I} + 1$.
\end{cor}

\noindent
{\it Proof.}\hspace{0.5ex}
This follows from \ref{3a} since
$(\,\ic{I}\,)^2 = Q\ic{I}$ by \cite[5.1]{Hn}
(or \cite[5.4]{LT}).

\begin{cor}\label{3c}
Let $\gp$ be a prime ideal of $A$
with $\height{}{\gp} = g \geq 2$.
Let $Q = (a_1, a_2, \dots, a_g)$ be an ideal
generated by a regular sequence contained in 
the $k$-th symbolic power $\gp^{(k)}$ of $\gp$
for some $k \geq 2$.
Then we have $\rn{Q}{I} \leq \nog{A}{(Q : \gp^{(k)}) / Q} + 1$
for any ideal $I$ with
$Q \subseteq I \subseteq Q : \gp^{(k)}$,
if one of the following three conditions holds {\rm ;}
{\rm (i)}
$A_\gp$ is not a regular local ring,
{\rm (ii)}
$A_\gp$ is a regular local ring and $g \geq 3$,
{\rm (iii)}
$A_\gp$ is a regular local ring,
$g = 2$, and $a_i \in \gp^{(k + 1)}$
for any $1 \leq i \leq g$.
\end{cor}

\noindent
{\it Proof.}\hspace{0.5ex}
This follows from \ref{3a} since
$(Q : \gp^{(k)})^2 = Q(Q : \gp^{(k)})$
by \cite[3.1]{W}.

\begin{cor}\label{3d}
Let $(A, \gm)$ be a Buchsbaum local ring.
Assume that the multiplicity of $A$
with respect to $\gm$ is $2$ and $\dep{}{A} > 0$.
Then, for any parameter ideal $Q$ in $A$
and an ideal $I$ with $Q \subseteq I \subseteq Q : \gm$,
we have $\rn{Q}{I} \leq \nog{A}{(Q : \gm) / Q} + 1$.
\end{cor}

\noindent
{\it Proof.}\hspace{0.5ex}
This follows from \ref{3a} since
$(Q : \gm)^2 = Q(Q : \gm)$ by \cite[1.1]{GS}.

\vspace{1.2em}
In order to state the last corollary,
let us recall the definition of Hilbert coefficients.
Let $(A, \gm)$ be a $d$-dimensional Noetherian local ring
and $I$ an $\gm$-primary ideal.
Then there exists a family
$\{\,\mult{i}{I}\,\}_{0 \leq i \leq d}$  of integers such that
\[
\length{A}{A / I^{n + 1}} = \sum_{i = 0}^d\,
(-1)^i\, \mult{i}{I}\, \binom{n + d - i}{d - i}
\]
for $n \gg 0$.
We call $\mult{i}{I}$ the $i$-th Hilbert coefficient of $I$.
On the other hand,
if $A$ is an analytically unramified local ring,
then $\{\,\ic{I^n}\,\}_{n \geq 0}$ is a Hilbert filtration
(cf. \cite{GR}), and so there exists a family
$\{\,\nmult{i}{I}\,\}_{0 \leq i \leq d}$ of integers such that
\[
\length{A}{A / {\ic{I^{n + 1}}}} = \sum_{i = 0}^d\,
(- 1)^i\, \nmult{i}{I}\, \binom{n + d - i}{d - i}
\]
for $n \gg 0$.
As is proved in \cite[1.5]{R},
if $A$ is a two-dimensional Cohen-Macaulay local ring,
then we have
\[
\rn{Q}{I} \leq \mult{1}{I} - \mult{0}{I} + \length{A}{A / I} + 1
\]
for any minimal reduction $Q$ of $I$.
We can generalize this result as follows.

\begin{cor}\label{3f}
Let $(A, \gm)$ be a two-dimensional Cohen-Macaulay local ring
with infinite residue field and $I$ an $\gm$-primary ideal
with a minimal reduction $Q$.
Then we have the following inequalities.
\begin{itemize}
\item[{\rm (1)}]
$\rn{Q}{I} \leq \mult{1}{J} - \mult{0}{J} + \length{A}{A / I} + 1$
for any ideal $J$ such that $I \subseteq J \subseteq \ic{I}$.
\item[{\rm (2)}]
$\rn{Q}{I} \leq \nmult{1}{I} - \nmult{0}{I} + \length{A}{A / I} + 1$,
if $A$ is analytically unramified.
\end{itemize}
\end{cor}

\noindent
{\it Proof.}\hspace{0.5ex}
(1)\hspace{0.5ex}
Setting $F_n = \rrc{J^n}$ for any $n \geq 0$ in \ref{1b}, we get
\begin{eqnarray*}
\rn{Q}{I} & \leq  & 
  1 + \nog{A}{\rrc{J} / I} + \sum_{n \geq 2}
                             \nog{A}{\rrc{J^n} / Q\rrc{J^{n - 1}}} \\
          & \leq  &
  1 + \length{A}{\rrc{J} / I} + \sum_{n \geq 2}
                                \length{A}{\rrc{J^n} / Q\rrc{J^{n - 1}}} \\
          & = & 
  \sum_{n \geq 1}\,\length{A}{\rrc{J^n} / Q\rrc{J^{n - 1}}} -
             \length{A}{I / Q} + 1\,.
\end{eqnarray*}
Because $\mult{1}{J} =
\sum_{n \geq 1}\,\length{A}{\rrc{J^n} / Q\rrc{J^{n - 1}}}$
by \cite[1.10]{GR} and
\[
\length{A}{I / Q} = \length{A}{A / Q} -
\length{A}{A / I} = \mult{0}{J} - \length{A}{A / I}\,,
\]
the required inequality follows.

(2)\hspace{0.5ex}
Similarly as the proof of (1),
setting $F_n = \ic{I^n}$ for any $n \geq 0$ in \ref{1b}, we get
\[
\rn{Q}{I} \leq \sum_{n \geq 1}\,\length{A}{\ic{I^n} / Q\ic{I^{n - 1}}}
- \length{A}{I / Q} + 1\,.
\]
Because the depth of the associated graded ring of the filtration
$\{\,\ic{I^n}\,\}_{n \geq 0}$ is positive,
we have $\nmult{1}{I} = \sum_{n \geq 1}\,\length{A}{\ic{I^n} / Q\ic{I^{n - 1}}}$
by \cite[1.9]{GR}.
Hence we get the required inequality as
$\length{A}{I / Q} = \nmult{0}{I} - \length{A}{A / I}$.

\section{Example}
In this section we give an example which shows that
the maximum value stated in \ref{3a} can be reached.
It provides an example in the case where
$\dim {A / I} > 0$.

\begin{ex}\label{4a}
Let $n \geq 3$ be an integer and
$S = k[X_0, X_1, \dots , X_n]$ be the polynomial ring
with $n + 1$ variables over a field $k$.
Let $A = S / \ga$,
where $\ga$ is the ideal of $S$ generated by
the maximal minors of the matrix
\[
\left(\begin{array}{llll}
X_0 & X_1 & \cdots & X_{n - 1} \\
X_1 & X_2 & \cdots & X_n
\end{array}\right)\,.
\]
We denote the image of $X_i$ in $A$ by $x_i$ for $0 \leq i \leq n$.
It is well known that $A$ is a two-dimensional Cohen-Macaulay graded ring
with the graded maximal ideal $\gm = (x_0, x_1, \dots , x_n)$.
\begin{itemize}
\item[{\rm (1)}]
Let $I = (x_0, x_1, x_n)$ and $Q = (x_0, x_n)$.
Then we have $\gm^2 = Q\gm$, $\nog{A}{\gm / I} = n - 2$, and
$\rn{Q}{I} = n - 1$.
\item[{\rm (2)}]
Let $I = (x_0, x_1, x_{n - 1})$, $J = (x_0, x_1, \dots , x_{n - 1})$,
and $Q = (x_0, x_{n - 1})$.
Then we have $\dim A / I = 1$, $J^2 = QJ$,
$\nog{A}{J / I} = n - 3$, and $\rn{Q}{I} = n - 2$.
\end{itemize}
\end{ex}

\noindent
{\it Proof.}\hspace{0.5ex} (1) \hspace{0.5ex}
Let $0 \leq i \leq j \leq n$.
If $i = 0$ or $j = n$,
then $x_ix_j \in Q\gm$.
On the other hand,
if $i > 0$ and $j < n$,
then the determinant of the matrix
\[
\left(\begin{array}{ll}
X_{i - 1} & X_{j} \\
X_i & X_{j + 1}
\end{array}\right)
\]
is contained in $\ga$, and so $x_ix_j = x_{i - 1}x_{j+1}$.
Hence we can show that $x_ix_j \in Q\gm$
for any $0 \leq i \leq j \leq n$ by descending induction on $j - i$.
Thus we get $\gm^2 = Q\gm$.
It is obvious that $\nog{A}{\gm / I} = n - 2$.
Therefore $I^n = QI^{n - 1}$ by \ref{3a}
(In fact, we have
${x_1}^n = {x_1}^{n - 2} \cdot {x_1}^2 = {x_1}^{n - 2} \cdot x_0x_2 =
x_0{x_1}^{n - 3} \cdot x_1x_2 = x_0{x_1}^{n - 3} \cdot x_0x_3 =
{x_0}^2{x_1}^{n - 4} \cdot x_1x_3 = \cdots =
{x_0}^{n - 2} \cdot x_1x_{n - 1} = {x_0}^{n - 2} \cdot x_0x_n =
{x_0}^{n - 1}x_n \in Q^n \subseteq QI^{n - 1}$).
In order to prove $\rn{Q}{I} = n - 1$,
we show ${x_1}^{n - 1} \not\in QI^{n - 2}$.
For that purpose we use the isomorphism
\[
\varphi : A \longrightarrow k[\,\{\,s^{n - i}t^i\,\}_{0 \leq i \leq n}\,]
\]
of $k$-algebras such that $\varphi(x_i) = s^{n - i}t^i$
for $0 \leq i \leq n$,
where $s$ and $t$ are indeterminates.
We have to show
${\varphi(x_1)}^{n - 1} \not\in \varphi(Q){\varphi(I)}^{n - 2}$.
Because $\varphi(I) = (s^n, s^{n - 1}t, t^n)$,
we get
\[
{\varphi(I)}^\ell \subseteq
(\, \{\, s^{\alpha n - \beta}t^{(\ell - \alpha)n + \beta} \, \mid \,
0 \leq \alpha \leq \ell \, , \, 0 \leq \beta \leq \alpha \, \} \, )
\]
for any $\ell \geq 1$ by induction on $\ell$, and so
\[
\varphi(Q){\varphi(I)}^{n - 2} \subseteq
( \, \{ \, s^{(\alpha + 1)n - \beta}t^{(n - 2 - \alpha)n + \beta} \, , \,
s^{\alpha n - \beta}t^{(n - 1 - \alpha)n + \beta} \, \mid \,
0 \leq \alpha \leq n - 2 \, , \, 0 \leq \beta \leq \alpha \, \} \, )\,.
\]
Therefore, if $\varphi(x_1)^{n - 1} = (s^{n - 1}t)^{n - 1} =
s^{(n - 1)^2}t^{n - 1} \in \varphi(Q){\varphi(I)}^{n - 2}$,
one of the following two cases
\begin{itemize}
\item[{(i)}]
$(\alpha + 1)n - \beta \leq (n - 1)^2$ and $(n - 2 - \alpha)n + \beta \leq n - 1$, or
\item[{(ii)}]
$\alpha n - \beta \leq (n - 1)^2$ and $(n - 1 - \alpha)n + \beta \leq n - 1$
\end{itemize}
must occur for some $\alpha$ and $\beta$ with
$0 \leq \alpha \leq n - 2$ and $0 \leq \beta \leq \alpha$.
Suppose that the case (i) occured.
Then we have
\[
(\alpha + 1)n - \beta \leq (n - 1)n - (n - 1) \hspace{1ex} \mbox{and} \hspace{1ex}
(n - 2 - \alpha)n \leq n - 1 - \beta\,.
\]
As the first inequality implies
\[
n - 1 - \beta \leq (n - 1)n - (\alpha + 1)n = (n - 2 - \alpha)n\,,
\]
it follows that
\[
n - 1 - \beta = (n - 1)n - (\alpha + 1)n\,,
\]
and so
\[
\alpha n - \beta = n^2 - 3n + 1\,.
\]
Then, as $\alpha n > n^2 - 3n = (n - 3)n$, we have
$n - 3 < \alpha \leq n - 2$, which implies $\alpha = n - 2$.
Thus we get
\[
(n - 2)n - \beta = n^2 - 3n + 1\,,
\]
and so $\beta = n - 1$, which contradicts to $\beta \leq \alpha$.
Therefore the case (ii) must occur.
Then we have
\[
\alpha n - \beta \leq (n - 1)n - (n - 1) \hspace{1ex} \mbox{and} \hspace{1ex}
(n - 1 - \alpha)n \leq n - 1 - \beta\,.
\]
As the first inequality implies
\[
n - 1 - \beta \leq (n - 1)n - \alpha n = (n - 1 - \alpha)n\,,
\]
it follows that
\[
n - 1 - \beta = (n - 1)n - \alpha n\,,
\]
and so
\[
\alpha n - \beta = n^2 - 2n + 1\,.
\]
Then, as $\alpha n > n^2 - 2n = (n - 2)n$,
we get $\alpha > n - 2$,
which contradicts to $\alpha \leq n - 2$.
Thus we have seen that ${x_1}^{n - 1} \not\in QI^{n - 2}$.

(2)\hspace{0.5ex}
Let $\gb = (X_0, X_1, \dots , X_{n - 1})S$.
Then $\ga \subseteq \gb$,
and so $\gb$ is the kernel of the canonical surjection
$S \longrightarrow A / J$.
Hence $A / J \cong k[\,X_n\,]$,
which implies $\dim {A / J} = 1$.
Let $0 \leq i \leq j \leq n - 1$.
If $i = 0$ or $j = n - 1$,
then $x_ix_j \in QJ$.
On the other hand,
if $i > 0$ and $j < n$, then $x_ix_j = x_{i - 1}x_{j + 1}$.
Hence we can show that $x_ix_j \in QJ$ for any
$0 \leq i \leq j \leq n - 1$ by descending induction on $j - i$.
Thus we get $J^2 = QJ$.
It is obvious that $\nog{A}{J / I} = n - 3$.
Therefore $I^{n - 1} = QI^{n - 2}$ by \ref{3a}.
This means $\dim {A / I} = \dim {A / Q} = \dim {A / J} = 1$.
In order to prove $\rn{Q}{I} = n - 2$,
we show ${x_1}^{n - 2} \not\in QI^{n - 3}$.
For that purpose we use again the isomorphism
$\varphi$ stated in the proof of (1).
Although we have to prove
${\varphi(x_1)}^{n - 2} \not\in \varphi(Q){\varphi(I)}^{n - 3}$,
it is enough to show
\[
(s^{n - 1}t)^{n - 2} \not\in
(s^n, st^{n - 1})(s^n, s^{n - 1}t, st^{n - 1})^{n - 3}B\,,
\]
where $B = k[s, t]$.
Because
\[
(s^{n - 1}t)^{n - 2} = s^{n - 2} \cdot (s^{n - 2}t)^{n - 2}
\]
in $B$ and
\[
(s^n, st^{n - 1})(s^n, s^{n - 1}t, st^{n - 1})^{n - 3}B =
s^{n - 2} \cdot
(s^{n - 1}, t^{n - 1})(s^{n - 1}, s^{n - 2}t, t^{n - 1})^{n - 3}B\,,
\]
we would like to show
\[
(s^{n - 2}t)^{n - 2} \not\in
(s^{n - 1}, t^{n - 1})(s^{n - 1}, s^{n - 2}t, t^{n -1})^{n - 3}B\,.
\]
However, it can be done by the same argument as the proof of
\[
(s^{n - 1}t)^{n - 1} \not\in
(s^n, t^n)(s^n, s^{n - 1}t, t^n)^{n - 1}\,,
\]
and hence we have proved (2).

\end{document}